\documentclass[12pt]{article}
\usepackage{amssymb}
\usepackage{url}
\newtheorem{theoreme}{Theorem}[section]
\newtheorem{defin}[theoreme]{Definition}
\newtheorem{remark}[theoreme]{Remark}
\newtheorem{prop}[theoreme]{Proposition}
\newtheorem{lemme}[theoreme]{Lemma}
\newtheorem{corol}[theoreme]{Corollary}

\newenvironment{demo}{\noindent {\sl Proof}. \ }{\qed}

\newtheorem{stheoreme}{Theorem}[subsection]
\newtheorem{sdefin}[stheoreme]{Definition}
\newtheorem{sremark}[stheoreme]{Remark}
\newtheorem{sprop}[stheoreme]{Proposition}
\newtheorem{slemme}[stheoreme]{Lemma}

\def\ssum{\sum\limits}

\font\tenCal=cmsy10
\newfam\Calfam
\textfont\Calfam=\tenCal

\def\qed{\hfill{$\sqcap\!\!\!\!\sqcup$}}

\def\t{\tilde}
\def\wt{\widetilde}

\def\calA{{\mathcal A}}
\def\derA{{\mathcal D}_{\sigma}({\mathcal A})}
\def\innA{{\mathcal I}{\rm nn}_{\sigma}({\mathcal A})}

\def\l{\lambda}
\def\a{\alpha}

\def\d{\delta}

\title{Quasi-Lie structure of $\sigma$-derivations of ${\mathbb C}[t^{\pm 1}]$}
\author{L. Richard\footnote{University of Edinburgh,
School of Mathematics, JCMB - King's Buildings, EH9 3JZ Edinburgh,
United Kingdom. lionel.richard@ed.ac.uk}, \ S.
Silvestrov\footnote{Centre for Mathematical Sciences, Lund
University, Box 118, SE-221 00 Lund, Sweden.
Sergei.Silvestrov@math.lth.se}}

\begin{document}

\maketitle

\begin{abstract}
Hartwig, Larsson and the second author in [J. Algebra, 295 (2006)] defined a
bracket on $\sigma$-derivations of a commutative algebra. We show
that this bracket preserves inner derivations, and based on this
obtain some structural results on $\sigma$-derivations on Laurent
polynomials in one variable. \\

\end{abstract}

\section{Introduction}

In \cite{HLS,LS1,LS2,LS3} a new class of algebras called quasi-Lie algebras
and its subclasses, quasi-hom-Lie algebras and hom-Lie algebras,
have been introduced. An important characteristic feature of those
algebras is that they obey some deformed or twisted versions of
skew-symmetry and Jacobi identity with respect to some possibly
deformed or twisted bilinear bracket multiplication. Quasi-Lie
algebras include color Lie algebras, and in particular Lie algebras
and Lie superalgebras, as well as various interesting quantum
deformations of Lie algebras, in particular of the Heisenberg Lie
algebra, oscillator algebras, $sl_2$ and other finite-dimensional
Lie algebras as well as of infinite-dimensional Lie algebras of Witt
and Virasoro type applied in physics within the string theory,
vertex operator models, quantum scattering, lattice models and other
contexts (see
\cite{ChaiIsKuLuk,ChaiElinPop,ChaiPopPres,ChaiKuLukPopPresn,
CurtrZachos1,DamKu,DelEtinQuantFieldStrings,DiFranMathiSen,FLM,
Fuchs1,Fuchs2,HLS,HelSil-book,Hu,Jakob,JakobLee,K} and references
therein). Many of these quantum deformations of Lie algebras can be
shown to play role of underlying algebraic objects for calculi of
twisted, discretized or deformed derivations and difference type
operators and thus in corresponding general non-commutive
differential calculi.

In \cite[Theorem 5]{HLS}, it was proved  that under some general
assumptions, when derivations are replaced by twisted derivations,
vector fields become replaced by twisted vector fields closed under
a natural twisted skew-symmetric bracket multiplication satisfying a
twisted 6-term Jacobi identity generalizing the usual Lie algebras
3-term Jacobi identity that is recovered when no twisting is present
(see Theorem \ref{thmbracket}). This theorem is shown to be
instrumental for construction of various examples and classes of
quasi-Lie algebras. Both known and new one-parameter and
multi-parameter deformations of Witt and Virasoro algebras and other
Lie and color Lie algebras has been constructed within this
framework in \cite{HLS,LS1,LS2,LS3}.

In this article, we gain further insight in the particular class of
quasi Lie algebra deformations of the Witt algebra, introduced in
\cite{HLS} via the general twisted bracket construction, and
associated with twisted discretization of derivations generalizing
the Jackson $q$-derivatives to the case of twistings by general
endomorphisms of Laurent polynomials. In Section
\ref{sec:generalfacts} we present necessary definitions, facts and
constructions on $\sigma$-derivations that are central for this
article. In Proposition \ref{prop132} we observe that inner
derivations are stable under the bracket defined in \cite{HLS}. In
Section \ref{secufd}, we present a characterization of the set of
inner derivations for UFD (Proposition \ref{innufd}), and also
general inclusions concerning sets of inner derivations and image
and pre-image subsets with respect to the twisted bracket
(Proposition \ref{prop:incusions}). In Section
\ref{sec:Non-linqdefWitt} we develop the preceding framework for a
particular important UFD, the algebra ${\mathcal A}={\mathbb
C}[t^{\pm 1}]$ of Laurent polynomials in one variable. For this
specialization more deep and precise results can be obtained. In
this case, $\derA$ with the twisted bracket is the deformation of
the Witt algebra within the class of quasi-hom-Lie algebras in the
sense of \cite{HLS}. In Theorem \ref{thmdecompo} we show that the
space of $\sigma$-derivations can be decomposed into a direct sum of
the space of inner $\sigma$-derivations and a dependent on $\sigma$
finite number of one-dimensional subspaces. In Theorem
\ref{propgrad}, we show that for arbitrary $\sigma$ the ${\mathbb
Z}$-gradation of this non-linearly $q$-deformed Witt algebra with
coefficients in ${\mathbb C}$ becomes a ${\mathbb Z}/d{\mathbb
Z}$-gradation with coefficients in ${\mathbb C}[T^{\pm 1}]$ for some
element $T$ of $A$. The usual $q$-deformed Witt algebra associated
to ordinary Jackson $q$-derivative, the automorphism case,
corresponds to $d=0$, when all $\sigma$-derivations are inner. In
Subsection \ref{parmodulo}, we provide more detailed description of
what relations for the bracket in the non-linearly $q$-deformed Witt
algebra become modulo inner $\sigma$-derivations. Finally, in
Subsection \ref{subsec:qWittS1}, we describe stabilizer-like subsets
in detail for the non-linearly $q$-deformed Witt algebra.

Throughout this article, ${\mathcal A}$ will denote an associative
and unital algebra over the field ${\mathbb C}$ of complex numbers.
The algebra ${\mathcal A}$ will often be  assumed to be commutative,
but we will sometimes mention more general results concerning
non-commutative algebras, so we may precise our assumptions on
${\mathcal A}$ every time.

\section{Some general facts on $\sigma$-derivations}
\label{sec:generalfacts}
\subsection{Definitions}
We recall here some basic definitions and facts concerning $\sigma$-derivations.
On this subject and more generally on Ore extensions one may see the reference book \cite{MCR}, or section 1.7 in \cite{Kabook}.
\begin{sdefin}
Let ${\mathcal A}$ be an algebra, and $\sigma$ an endomorphism of ${\mathcal A}$.
A $\sigma$-derivation is a linear map $D$ satisfying $D(ab)=\sigma(a)D(b)+D(a)b$
for all $a,b\in {\mathcal A}$. We denote the set of all $\sigma$-derivations by ${\mathcal D}_{\sigma}({\mathcal A})$.
\end{sdefin}
{\bf Example.} It is easy to check that for any $p\in{\mathcal A}$,
the ${\mathbb C}$-linear map $\Delta_p$ defined by $\Delta_p(a)=pa-\sigma(a)p$
for all $a\in{\mathcal A}$ is a $\sigma$-derivation of ${\mathcal A}$.
Note that if ${\mathcal A}$ is commutative, then we have $\Delta_p=p({\rm id}-\sigma)$.
\begin{sdefin}\label{definn}
The map $\Delta_p$ defined above is called the {\rm inner} $\sigma$-de\-ri\-va\-tion associated to $p$.
The set of all inner derivations of ${\mathcal A}$ will be denoted ${\mathcal I}nn_{\sigma}({\mathcal A})$.
\end{sdefin}
For any map $\tau:\calA \rightarrow \calA$ denote ${\rm
Ann}(\tau)=\{a\in \calA\ |\ a\tau(b)=0\ \forall\ b\in\calA\}$, the
left annihilator ideal of $\tau$. In particular if ${\mathcal A}$ is
commutative then this is a two sided ideal, and also
$\Delta_p=\Delta_q\iff(p-q)\in{\rm Ann}({\rm id}-\sigma)$.
\medskip

The $\sigma$-derivations play a crucial role in the definition of Ore extensions, that we recall here.

\begin{sdefin}
Let ${\mathcal A}$ be an algebra, $\sigma$ an endomorphism of
${\mathcal A}$, and $\Delta$ a $\sigma$-derivation. Then the Ore
extension $R={\mathcal A}[X;\sigma, \Delta]$ is the algebra such
that:
\begin{itemize}
\item ${\mathcal A}$ is a sub-algebra of $R$;
\item $R$ is a free ${\mathcal A}$-module with basis $\{X^n,\ n\in{\mathbb N}\}$;
\item the multiplication is defined in $R$ by the rule $Xa=\sigma(a)X+\Delta(a)$ for all $a\in {\mathcal A}$.
\end{itemize}
\end{sdefin}

The  following facts can be easily checked; they appear in
\cite[Lemmas 1.5 and 2.4.]{G}.
\begin{slemme}\label{goodearl}
Let $\sigma$ be an endomorphism of an algebra ${\mathcal A}$ and
$\Delta_p$ be an inner derivation corresponding to an element $p\in
A$.
\begin{enumerate}
\item Set $p\in\calA$, and $\Delta_p\in\innA$.
Then the identity map on ${\mathcal A}$ extends to an isomorphism
$\tau$ between the Ore extensions ${\mathcal A}[X;\sigma, \Delta_p]$
and ${\mathcal A}[Y;\sigma]$, defined by $\tau(X)=Y+p$ sending $X$
to $Y+p$.
\item Assume ${\mathcal A}$ is commutative. Then for all $\Delta\in\derA$, and for all $a,b\in{\mathcal A}$ one has
$$(a-\sigma(a))\Delta(b)=(b-\sigma(b))\Delta(a).$$
\end{enumerate}
\end{slemme}

The first statement of this lemma provides one of the reasons why we
get interested in $\sigma$-derivations up to inner in the following
sections.

\subsection{A bracket on $\sigma$-derivations}

From now on ${\mathcal A}$ is supposed to be a {\it commutative}
algebra. Then ${\mathcal D}_{\sigma}({\mathcal A})$ becomes a left
${\mathcal A}$-module by $(a\Delta)(r)=a\Delta(r)$ for all $a,r\in
{\mathcal A}$. Note that in the non-commutative case this operation
makes ${\mathcal D}_{\sigma}({\mathcal A})$ a left module only over
the center of ${\mathcal A}$.

Now we fix a $\sigma$-derivation $\Delta$, and consider the cyclic
left ${\mathcal A}$-submodule of ${\mathcal D}_{\sigma}({\mathcal
A})$ generated by $\Delta$. The interest of considering cyclic
sub-modules is reinforced by the following result proved by Hartwig,
Larsson and the second author in \cite{HLS}, Theorem 2, in the case
where ${\mathcal A}$ is a unique factorization domain.
\begin{stheoreme}[\cite{HLS}]\label{ufd}
Let $\sigma$ be an endomorphism of a unique factorization domain ${\mathcal A}$, and $\sigma\neq{\rm id}$.
Then ${\mathcal D}_{\sigma}({\mathcal A})$ is free of rank one as an ${\mathcal A}$-module with generator
$$\Delta={{{\rm id}-\sigma}\over{g}}, {\rm with} \ \ g={\rm gcd}(({\rm id}-\sigma)({\mathcal A})).$$
\end{stheoreme}

They also define in \cite{HLS}, Theorem 5, a bracket on this  cyclic
sub-module and prove the following results.
\begin{stheoreme}[\cite{HLS}]\label{thmbracket}
Let $\sigma$ be an endomorphism of a commutative algebra ${\mathcal A}$, and $\sigma\neq{\rm id}$.
Let $\Delta\in{\mathcal D}_{\sigma}({\mathcal A})$ be a $\sigma$-derivation such that :
\begin{itemize}
\item $\sigma({\rm Ann}(\Delta))\subseteq{\rm Ann}(\Delta)$;
\item $\exists\ \delta\in{\mathcal A}$ such that $\Delta\circ\sigma=\delta \sigma \circ \Delta$.
\end{itemize}
then the map
$$[\cdot,\cdot]_{\sigma}:\calA\Delta\times\calA\Delta \to
    \calA\Delta$$
defined by setting
\begin{equation} \label{eq:GenWittProdDef}
    [a\Delta,b\Delta]_{\sigma}=(\sigma(a)\Delta)\circ(b\Delta)-(\sigma(b)\Delta)
    \circ(a\Delta), \quad\textrm{for }a,b\in\calA,
\end{equation}
where $\circ$ denotes composition of functions, is a well-defined
${\mathbb C}$-algebra product on the ${\mathbb C}$-linear space $\calA\Delta$,
satisfying the following identities \mbox{for $a,b,c\in\calA$:}
\begin{equation}\label{eq:GenWittProdFormula}
    [a\Delta, b\Delta]_\sigma=\big(\sigma(a)\Delta(b)-\sigma(b)\Delta(a)\big)\Delta,
\end{equation}
\begin{equation}\label{eq:GenWittSkew}
    [a\Delta, b\Delta]_\sigma=-[b\Delta, a\Delta]_\sigma.
\end{equation} In addition,
\begin{equation} \label{eq:GenWittJacobi}
\begin{array}{c}
   [\sigma(a)\Delta,[b\Delta,c\Delta]_\sigma]_\sigma+\delta[a\Delta,[b\Delta,c\Delta]_\sigma]_\sigma+\\
   +[\sigma(b)\Delta,[c\Delta,a\Delta]_\sigma]_\sigma+\delta[b\Delta,[c\Delta,a\Delta]_\sigma]_\sigma+\\
   +[\sigma(c)\Delta,[a\Delta,b\Delta]_\sigma]_\sigma+\d[c\Delta,[a\Delta,b\Delta]_\sigma]_\sigma=0.
\end{array}
\end{equation}
\end{stheoreme}
%\qed

\begin{sremark} {\rm \begin{enumerate}
\item The second condition, under the extra assumption that $\delta=q\in {\mathbb C}^*$, is the definition of $q$-skew derivations given in \cite{G}. These particular $\sigma$-derivations play an important role in  quantum groups, see for instance \cite{C} or the reference book \cite{BG} and references therein.
Note that in this case Formula (\ref{eq:GenWittJacobi}) can be written as a 3-term Jacobi-like identity.
\item     The identity (\ref{eq:GenWittProdFormula}) is just a formula
expressing the
product defined in
    (\ref{eq:GenWittProdDef}) as an element of $\calA\Delta$. Identities
(\ref{eq:GenWittSkew}) and
    (\ref{eq:GenWittJacobi}) are  expressing, respectively,
skew-symmetry and a generalized ($(\sigma,\delta)$-twisted) Jacobi identity
for the product defined by (\ref{eq:GenWittProdDef}).
\end{enumerate}}
\end{sremark}

\subsection{Inner $\sigma$-derivations}
We recall from Definition \ref{definn} that a $\sigma$-derivation
$\wt\Delta$ is inner if and only if there exists an element
$p\in{\mathcal A}$ such that $\wt\Delta(a)=pa-\sigma(a)p$ for all
$a\in{\mathcal A}$. Because ${\mathcal A}$ is commutative, it is
easy to see that if $\Delta$ is inner, then $a\Delta$ is inner for
all $a\in{\mathcal A}$, so that ${\mathcal
A}\Delta\subseteq{\mathcal I}nn_{\sigma}({\mathcal A})$. So we
mainly get interested in the case where $\Delta$ itself is not
inner.

\smallskip

Now we  prove that inner derivations are stable under the bracket $[.,.]_{\sigma}$.

\begin{sprop} \label{prop132}
Let $\sigma$ be an endomorphism of a commutative algebra ${\mathcal
A}$, and $\Delta\in{\mathcal D}_{\sigma}({\mathcal A})$. Set
$a,b\in{\mathcal A}$ such that $a\Delta=p({\rm id}-\sigma)$ and
$b\Delta=q({\rm id}-\sigma)$ are inner. Then
$[a\Delta,b\Delta]_{\sigma}$ is inner. More precisely, we have:
$[a\Delta,b\Delta]_{\sigma}=c({\rm id}-\sigma)$, with
$c=\Delta(b)p-\Delta(a)q$, and so
$[a\Delta,b\Delta]_{\sigma}=(\Delta(b)a-\Delta(a)b)\Delta$.
\end{sprop}
\begin{demo}
For all $r\in {\mathcal A}$ one has: $$[a\Delta,b\Delta]_{\sigma}(r)=(\sigma(a)\Delta(b)-\sigma(b)\Delta(a))\Delta(r)=\sigma(a)\Delta(b)\Delta(r)-\sigma(b)\Delta(a)\Delta(r).$$
From Lemma \ref{goodearl} we have $(r-\sigma(r))\Delta(a)=(a-\sigma(a))\Delta(r)$.
Using the fact that $a\Delta(r)=p(r-\sigma(r))$
and $\calA$ is commutative we get then $\sigma(a)\Delta(r)=(r-\sigma(r))(p-\Delta(a))$.
In the same way one proves that $\sigma(b)\Delta(r)=(r-\sigma(r))(q-\Delta(b))$. So:
$$\begin{array}{rcl}
[a\Delta,b\Delta]_{\sigma}(r)&=&\Delta(b)(r-\sigma(r))(p-\Delta(a))-\Delta(a)(r-\sigma(r))(q-\Delta(b))\\
&=&(\Delta(b)p-\Delta(a)q)(r-\sigma(r)).\end{array}$$
\end{demo}

%\medskip

\begin{sremark} \label{rkk}
{\rm \begin{enumerate} \item In the UFD case, since $a=gp$ and $b=gq$, then one gets $c=\sigma(g)(\Delta(q)p-\Delta(p)q)$.
\item $\innA$ appears as a sub-algebra of $(\derA,[.,.]_{\sigma})$. It is the sub-algebra considered in paragraph 3.2.1 of \cite{HLS}.
\item In the non-commutative case, under the assumption that $\Delta$ sends the center of ${\mathcal A}$ to itself,  we can prove in the same way  for $a,b\in Z({\mathcal A})$ and with the notations of
Definition \ref{definn}, that if $a\Delta=\Delta_p$ and $b\Delta=\Delta_q$ then $[a\Delta,b\Delta]_{\sigma}=\Delta_t$, with $t=\Delta(b)p-\Delta(a)q$.
\end{enumerate}}
\end{sremark}

\section{The UFD case.}\label{secufd}
We present in this section some general statements concerning
$\sigma$-derivations in the  UFD case.  In the next section we will
give some more precise and deep results in the particular case
${\mathcal A}={\mathbb C}[t^{\pm 1}]$. One can first precise Theorem
\ref{ufd} in the following way.
\begin{prop}\label{innufd}
 Let $\sigma$ be an endomorphism of ${\mathcal A}$  a unique factorization domain, $\sigma\neq{\rm id}$. Set $g={\rm gcd}(({\rm id}-\sigma)({\mathcal A}))$. Then
\begin{enumerate}
\item ${\mathcal D}_{\sigma}({\mathcal A})={\mathcal A}\Delta$, with $\Delta={{{\rm id}-\sigma}\over{g}}$;
\item the  $\sigma$-derivation $a\Delta$ is inner if and only if $g$ divides $a$.
In other words, $\innA=g{\mathcal A}\Delta$.
\end{enumerate}
\end{prop}
\begin{demo}
The first point is just Theorem \ref{ufd}. Then any
$\sigma$-derivation can be written $\widetilde\Delta=a\Delta$, with
$a\in{\mathcal A}$. Obviously if $a=bg$ then $a\Delta=b({\rm
id}-\sigma)$ is inner.

Conversely, assume that $a\Delta$ is inner. Then there is an element $b\in{\mathcal A}$ such that $a\Delta(r)=b(r-\sigma(r))$ for all $r\in{\mathcal A}$. Multiplying by $g$ one obtains $ag\Delta(r)=b(r-\sigma(r))$, i.e. $a(r-\sigma(r))=bg(r-\sigma(r))$ by definition of $\Delta$. Now choose $r$ such that $r-\sigma(r)\neq 0$, and conclude using the fact that ${\mathcal A}$ is a domain.
\end{demo}

\medskip

In particular, $\Delta$ itself is inner if and only if $g$ is a unit.

\medskip

In order to use the bracket $[.,.]_{\sigma}$ to understand what is ``between" $\derA$ and $\innA$, we define now some sub-spaces of $\derA$,  which will be  more precisely described in the next section for ${\mathcal A}={\mathbb C}[t^{\pm 1}]$.
Recall that $\derA={\mathcal A}\Delta$ and $\innA=g{\mathcal A}\Delta$.

\begin{defin}
$S^1={\rm Span}_{{\mathbb C}}[\innA,\innA]_{\sigma}$.
\end{defin}

\begin{remark}
{\rm \begin{enumerate} \item $S^1$ would be the usual derived Lie
algebra of $\innA$ for $\sigma={\rm id}$.
\item It follows from Proposition \ref{prop132} that one always have $S^1\subseteq \innA$.
\end{enumerate}}
\end{remark}

%\medskip
\begin{lemme} \label{s1sup}
$S^1\subseteq\sigma(g)g\calA\Delta$.
\end{lemme}
\begin{demo}
This inclusion relies on the remark following Proposition
\ref{prop132}. We are in the UFD case, and if $a,b\in \calA$ such
that $a\Delta,b\Delta\in\innA$, then $a=gp$ and $b=gq$ for some
$p,q\in\calA$, and $[a\Delta,b\Delta]_{\sigma}=\sigma(g)
(\Delta(q)p-\Delta(p)q)g\Delta$.  So
$[\innA,\innA]_{\sigma}\subseteq\sigma(g)g\calA\Delta$.
\end{demo}

\medskip

The following corollary is a direct consequence of the preceding Lemma and part 2 of Proposition \ref{innufd}.
\begin{corol}\label{cors1}
If $\sigma(g)$ is not a unit in ${\mathcal A}$, then $S^1\subset\innA$.
\end{corol}

%\medskip

Now we define two subspaces of $\derA$. Note that these definitions can be given for any algebra, no matter if it is a UFD or not.

\begin{defin}
$\wt S_1=\{\wt\Delta\in\derA |\ [\wt\Delta,S^1]\subseteq \innA\}$.
\end{defin}

Because $S^1\subseteq\innA$ it follows from Proposition \ref{prop132} that $\innA\subseteq \wt S_1$.
As we will see  in the case of $q$-deformed Witt algebras $\wt S_1$ is the whole space of $\sigma$-derivations,
so we define now a ``smaller'' space.

\begin{defin}
$S_1=\{\wt\Delta\in\derA\ |\ [\wt\Delta;S^1]_{\sigma}\subseteq S^1\}$.
\end{defin}
Because $S^1\subseteq \innA$ it is clear that $S_1\subseteq \wt S_1$.
This is the stabilizer of $S^1$ in $\derA$ with respect to
$[.,.]_\sigma$.

We can summerize the inclusions as follows.
\begin{prop} \label{prop:incusions}
Let ${\mathcal A}$ be a UFD. Then $S^1\subseteq\innA\subseteq
S_1\subseteq \wt S_1.$
\end{prop}

\section{Non-linearly $q$-deformed Witt algebras}
\label{sec:Non-linqdefWitt} We develop now the preceding framework
for a particular algebra ${\mathcal A}$, namely ${\mathbb C}[t^{\pm
1}]$, in order to obtain some more deep and precise results. We
study the case $\calA={\mathbb C}[t^{\pm 1}]$, so $\derA$ with the
twisted bracket is the deformation of the Witt algebra in the sense
of \cite{HLS}. The deformations of the Witt algebra are of
importance in mathematical physics (see \cite{HLS} and references
therein, and the introduction of the present work). Note that as
${\mathbb C}[t^{\pm 1}]$ is a UFD Proposition \ref{innufd} and
Proposition \ref{prop:incusions} apply.

Most of the articles on $\sigma$-derivations of ${\mathbb C}[t^{\pm
1}]$ are concerned with the case where $\sigma$ is an automorphism
(see for instance  \cite{K}, \cite{RS}). We don't assume here that
$\sigma$ is an automorphism, so it involves some power $s$ of $t$,
which as we will see plays a crucial role in the study of $\derA$.

\subsection{Some notations}

{\bf Degree, valuation.} For a Laurent polynomial $f(t)=\sum_{n=n_0}^{n_1}\alpha_nt^n$ with $\alpha_n\in{\mathbb C},\ \alpha_{n_0}\neq 0,\ \alpha_{n_1}\neq0$ we denote $\nu(f)=n_0$ the valuation of $f$ and ${\rm deg}(f)=n_1$ its degree.

\medskip

\noindent{\bf The endomorphism $\sigma$.} Because an endomorphism of algebra  sends units on units, the image of $t$ by $\sigma$ must be a monomial. So denote $\sigma(t)=qt^s$, with $q\in{\mathbb C}^*$ and $s\in{\mathbb Z}$. Note that $\sigma$ is injective if and only if $s\neq 0$, and surjective if and only if $s=1$ or $s=-1$.

\medskip

\noindent{\bf Generators of $\derA$.} If $\sigma\neq{\rm id}$ then by Theorem \ref{ufd} one has  $\derA=\calA\Delta$ with $\Delta(f)=(f-\sigma(f))/g$, $g={\textrm{gcd}}(({\textrm{id}}-\sigma)(\calA))$.
Then one can check (see \cite{HLS}, example 3.2) that $g=\a^{-1} t^{k-1}(t-qt^s)$ with $\a\in{\mathbb C}^*$ and $k\in{\mathbb Z}$. Since $g$ is defined up to a unit, then $\a$ and $k$ are arbitrary.
If $s\geq 1$ then choose $k=0$ and $\a=1$, so that $g(t)=1-qt^{s-1}$. If $s\leq 0$ then choose $k=-s+1$ and $\a=-q$, so that $g(t)=-q^{-1}(t^{1-s}-q)=1-q^{-1}t^{1-s}$. In both cases $g(t)$ is a (non-Laurent) polynomial of degree $|s-1|$ such that $g(0)=1$. So we can assume  without loss of generality that
$$g(t)=1-\l t^d,$$
 with $[\l=q, d=s-1]$ if $s\geq 1$ and $[\l=q^{-1},d=1-s]$ if $s<1$. So for $s\neq 1$ one has $d\geq1$ and $\l\neq 0$.
Note that with our conventions, for $\sigma={\rm Id}$ (i.e. $s=1=q$), one gets $g=0$.
We consider the following linear basis of $\derA$ (see Example 3.2 in \cite{HLS}): $d_n=-t^n\Delta$ for all ${\mathbb  Z}$.

The monomial $T=q t^{s-1}$ will play a crucial rule in the following. Note that  $g=1-T$ if $s\geq 1$, and $g=1-T^{-1}$ if $s<1$. Note also that  $\sigma(T)=T^s$, and that $T$ ``acts'' on $\derA$ in the following way: $Td_n=q d_{n+s-1}$.
At last, since  $\Delta(t)\neq 0$ and $\calA$ is a domain, we have ${\rm Ann}(\Delta)=\{0\}$.

\subsection{Decomposition of ${\mathcal D}_{\sigma}({\mathbb C}[t^{\pm1}])$}\label{pardecompo}

If $s=1$ then two cases may occur. First if $q=1$ then $\sigma$ is
the identity map, so one gets the usual Witt algebra, and we won't
consider this case here. Second case: if $q\neq 1$ then $g=1-q$ is
a unit, and Proposition \ref{innufd} implies that all
$\sigma$-derivations are inner, and $\derA=\innA=S^1=\wt S_1=S_1$.
As we get interested in the study of what happens ``between''
$\derA$ and $\innA$,  we shall assume that $s\neq 1$. Note that then
$g$ is not a unit in $\calA$, so thanks to Corollary \ref{cors1} we
have $S^1\subset \innA$.
%The situation is the following : we consider the set ${\mathcal D}_{\sigma}({\mathcal A})$
%of $\sigma$-derivations of the Laurent polynomials algebra ${\mathcal A}={\mathbb C}[t,t^{-1}]$,
%where $\sigma$ is the endomorphism defined by $\sigma(t)=qt^s$, $s\neq1$.

The vector space ${\mathcal D}_{\sigma}({\mathcal A})$ is made a
non-associative algebra thanks to the bracket $[.,.]_{\sigma}$
defined in Theorem \ref{thmbracket}. Since $g$ is not a unit, the
set ${\mathcal I}nn_{\sigma}({\mathcal A})$ of inner
$\sigma$-derivations is a proper subalgebra of $({\mathcal
D}_{\sigma}({\mathcal A}),[.,.]_{\sigma})$. Moreover, thanks to
Proposition \ref{innufd} we know that a $\sigma$-derivation
$f\Delta$ is inner if and only if $g$ divides $f$, that is
$\innA=g\calA\Delta$. This leads to the following result.

\begin{stheoreme} \label{thmdecompo}
Assume that $s\neq 1$. Then with the notations above
$$\derA={\mathbb C}d_0\oplus{\mathbb C}d_1\ldots\oplus{\mathbb C}d_{d-1}\oplus\innA.$$
\end{stheoreme}
\begin{demo}
Note that for any $n\in{\mathbb N}$ one has ${\mathbb C}d_n={\mathbb C}t^n\Delta$.
We first show that $\derA={\mathbb C}\Delta+{\mathbb C}t\Delta\ldots+{\mathbb C}t^{d-1}\Delta+\innA$.
Take any $\wt\Delta=f(t)\Delta\in\derA$, with $f(t)=\sum_{n=n_0}^{n_1}\alpha_nt^n$, with $\nu(f)=n_0$ and ${\rm deg}(f)=n_1$. Up to an inner derivation we can always assume that $\nu(f)\geq 0$. If not, consider $f_1=f-\a_{n_0}t^{n_0}g$ : we have $f\Delta=f_1\Delta +\a_{n_0}t^{n_0}g\Delta$, $\a_{n_0}t^{n_0}g\Delta\in\innA$ and $\nu(f_1)>\nu(f)$. If $\nu(f_1)\geq 0$ we are done, else we repeat this operation with $f_1$. Then after at most $\nu(f)$ iterations we have a polynomial $\t f$ such that $\nu(\t f)\geq 0$ and $g$ divides $f-\t f$, that is $(f\Delta-\t f\Delta)\in\innA$. So assume that $f\in {\mathbb C}[t]$.
Now we can make the usual Euclidian division of $f$ by  $g$ in ${\mathbb C}[t]$, and we obtain $f=q(t)g(t)+r(t)$, with ${\rm deg}(r)<{\rm deg}(g)=d$. Since $\innA=g\calA\Delta$ we are done.

Now we prove that this sum is a direct sum. Set $\a_0,\ldots,\a_{d-1}$ in ${\mathbb C}$ and $f\in{\mathbb C}[t^{\pm 1}]$ such that $\sum_{i=0}^{d-1} \a_it^i\Delta+fg\Delta=0$. Since ${\rm Ann}(\Delta)=\{ 0\}$ this implies $\sum_{i=0}^{d-1} \a_it^i+fg=0$. First we prove  that $f$ must be a non-Laurent polynomial, that is $\nu(f)\geq 0$. Assume on the contrary that $\nu(f)=n_0<0$, and $f$ has lowest degree term $\beta_{n_0}t^{n_0}\neq 0$. Then because $g=1-\l t^d$ with $d>0$, the term of lowest degree of $\sum_{i=0}^{d-1} \a_it^i+fg$ is $\beta_{n_0}t^{n_0}$, a contradiction since  $\sum_{i=0}^{d-1} \a_it^i+fg=0$.

Now the latest equality is nothing else but the Euclidian division of the $0$ polynomial by $g$ in ${\mathbb C}[t]$. By uniqueness we have $\a_i=0$ for all $i$ and $f=0$.
\end{demo}

\medskip

In Subsection \ref{parmodulo}  we will give a description of the
bracket in $\derA$ in terms of this decomposition, thanks to the
brackets computed in \cite{HLS}. But first we re-interpret these in
term of the element $T=qt^{s-1}$ defined at the beginning of this
section, and show that the algebra $\derA$ is graded by a finite
cyclic group.

\subsection{Grading of $\derA$}

We recall first the following from \cite[Theorem 8]{HLS}.

\begin{stheoreme} \label{th:qWittnonlinear} Let $\calA={\mathbb C}[t^{\pm 1}]$, and equip the
${\mathbb C}$-linear space
    $\derA = \bigoplus_{n\in{\mathbb Z}} {\mathbb C} d_n$
    with the bracket product
    $$[\,\cdot,\,\cdot]_\sigma\,:\, \derA\times \derA\longrightarrow
    \derA$$ defined on generators by (\ref{eq:GenWittProdDef}) as
    \begin{equation}\label{bracketdef1}
    [d_n,d_m]_\sigma=q^nd_{ns}d_m-q^md_{ms}d_n.
    \end{equation}
    This bracket satisfies
    defining commutation relations
$$        \begin{array}{l}
{[d_n,d_m]_\sigma}=\alpha{\rm{sign}}(n-m)\ssum_{l={\rm min}(n,m)}^{{\rm{max}}(n,m)-1}q^{n+m-1-l}
            d_{s(n+m-1)-(k-1)-l(s-1)}\\
            \qquad\textrm{for } n,m\geq 0;\\
            {[d_n,d_m]}_\sigma=\alpha\Big
           (\ssum_{l=0}^{-m-1}q^{n+m+l}d_{(m+l)(s-1)+ns+m-k}\\
       \qquad\qquad\qquad\qquad\qquad +\ssum_{l=0}^{n-1}q^{m+l}d_{(s-1)l+n+ms-k}\Big)\\
           \qquad\textrm{for } n\geq 0, m<0;
           \end{array}
$$
$$
            \begin{array}{l}
            {[d_n,d_m]}_\sigma=-\alpha\Big (\ssum_{l_1=0}^{m-1}q^{n+l_1}
        d_{(s-1)l_1+m+ns-k}\\
       \qquad\qquad\qquad\qquad\qquad  +\ssum_{l_2=0}^{-n-1}q^{m+n+l_2}d_{(n+l_2)(s-1)+n+ms-k}\Big )\\
            \qquad\textrm{for } m\geq 0, n<0;\\
           {[d_n,d_m]}_\sigma=\alpha{\rm{sign}}(n-m)\ssum_{l={\rm min}(-n,-m)}^{{\rm max}(-n,-m)-1}q^{n+m+l}d_{(m+n)s+(s-1)l-k}\\
            \qquad\textrm{for } n,m<0.
        \end{array}$$
    Furthermore, this bracket satisfies skew-symmetry
    $[d_n,d_m]_\sigma=-[d_m,d_n]_\sigma$ and a twisted Jacobi identity.
\end{stheoreme}

%\subsection{Graduation and relations modulo inner $\sigma$-derivations in case $s>1$} \label{sec42}
Note that for $s=1$ formula (\ref{bracketdef1}) gives the usual formula for the $q$-Witt algebra,
or the Witt algebra if $q=1$, that is
$$[d_n,d_m]_q=q^nd_nd_m-q^md_md_n.$$
Now for arbitrary $s\in{\mathbb Z}$, recall that $T=qt^{s-1}$, and
$\sigma(T)=T^s$. Assume $\sigma$ is not the identity map (so that
$T\neq 1$), and for all $n\in{\mathbb Z}$ note the $T$-integer
$\{n\}_T=\frac{T^n-1}{T-1}$. This is just a geometric sum, more
precisely one has $\{0\}_T=0$, $\{n\}_T=\sum_{k=0}^{n-1}T^k$ for
$n>0$ and $\{n\}_T=-\sum_{k=n}^{-1}T^k$ for $n<0$. Thanks to these
notations we shall rewrite the preceding formulas in the following
way.

\begin{sprop} \label{331}
For all $n<m\in{\mathbb Z}$ one has
$$[d_n;d_m]_{\sigma}=T^nd_nd_m-T^md_md_n;$$
\begin{equation}\label{bracketsgt1}
[d_n;d_m]_{\sigma}=-(T^n\sum_{k=0}^{m-n-1}T^k)d_{n+m}=(\{n\}_{T}-\{m\}_T)d_{n+m}.
\end{equation}
\end{sprop}
\begin{demo}
The first equation comes directly from (\ref{bracketdef1}) and the definitions of $T$ and of the $d_n$'s.
For the second one, we begin with the case $0\leq n<m$ :
$$[d_n,d_m]_{\sigma}=\sum_{l=n}^{m-1}q^{n+m-1-l}t^{s(n+m-1)+1-l(s-1)}\Delta;$$
which becomes after reindexing the summation with $k=m-1-l$ :
$$\begin{array}{rcl}[t^n\Delta,t^m\Delta]_{\sigma}&=&\ssum_{k=0}^{m-n-1}q^{n+k}t^{s(n+m-1)+1-(s-1)(m-1-k)}\Delta\\
&=&\ssum_{k=0}^{m-n-1}q^{n+k}t^{sn+m+k(s-1)}\Delta\\
&=&(q^nt^{sn+m}\ssum_{k=0}^{m-n-1}q^kt^{(s-1)k})\Delta\\
&=&t^{n+m}(qt^{s-1})^n\ssum_{k=0}^{m-n-1}(qt^{s-1})^k\Delta.\end{array}$$
The two other cases are treated in the same way, thanks to the following formulas :
if $n<0\leq m$ then
$$\begin{array}{rcl}[t^n\Delta,t^m\Delta]_{\sigma}&=&\ssum_{l=0}^{m-1}q^{n+l}t^{l(s-1)+m+ns}\Delta+\ssum_{l=0}^{-n-1}q^{m+n+l}t^{(n+l)(s-1)+n+ms}\Delta\\
&=&t^{n+m}(qt^{s-1})^n\ssum_{l=0}^{m-1}q^lt^{l(s-1)}\Delta\\
 & &      \qquad\qquad\qquad\qquad\qquad + \ssum_{l=m}^{m-n-1}q^{n+l}t^{(n+l-m)(s-1)+n+ms}\Delta\\
&=&t^{n+m}(qt^{s-1})^n\ssum_{k=0}^{m-n-1}(qt^{s-1})^k\Delta;\end{array}$$
and finally if $n<m<0$ then
$$\begin{array}{rcl}[t^n\Delta,t^m\Delta]_{\sigma}&=&\ssum_{l=-m}^{-n-1}q^{n+m+l}t^{(m+n)s+(s-1)l}\Delta\\
&=&\ssum_{l=0}^{m-n-1}q^{n+l}t^{(m+n)s+(s-1)(l-m)}\Delta\\
 &=&t^{n+m}(qt^{s-1})^n\ssum_{k=0}^{m-n-1}(qt^{s-1})^k\Delta.\end{array}$$
\end{demo}

\medskip

Let us remark here that the second expression in formula (\ref{bracketsgt1})
shows that these non-linearly deformed Witt algebras, constructed {\sl a priori}
in \cite{HLS} by taking any endomorphism of ${\mathbb C}[t^{\pm 1}]$ instead of the automorphism $t\to qt$, really ``look like'' the $q$-Witt algebra. More precisely, if one takes for $\sigma$ the automorphism $t\to qt$, then $s=1$, so  $d=0$, and $T=q$. Then the $T$-integers are the usual $q$-integers, and (\ref{bracketsgt1}) is the usual bracket of the $q$-Witt algebra, as defined for instance in \cite{K}, and leading for $q=1$ to the classical Witt algebra.
%\medskip

%\subsection{A gradation on $\derA$}

We show now that for arbitrary $\sigma$ the ${\mathbb Z}$-gradation
of the $q$-Witt algebra with coefficients in ${\mathbb C}$ becomes
a ${\mathbb Z}/d{\mathbb Z}$-gradation with coefficients in ${\mathbb C}[T^{\pm 1}]$.
Once again this is relevant, because the usual $q$-Witt algebra corresponds to the case where $d=0$.

\begin{stheoreme}\label{propgrad}
Let $\sigma$ be the endomorphism of ${\mathcal A}={\mathbb C}[t^{\pm
1}]$ defined by $\sigma(t)=qt^s$, with $q\in{\mathbb C}^*$ and
$s\in{\mathbb Z}$. Recall that $\derA$ is the space of
$\sigma$-derivations of ${\mathcal A}$, endowed with the bracket
$[.,.]_{\sigma}$ defined in Theorem {\rm \ref{thmbracket}}. Define
$d=|s-1|$, and note ${\mathbb Z}_d={\mathbb Z}/d{\mathbb Z}$, in
particular for $s=1$ one has ${\mathbb Z}/0={\mathbb Z}$. For any
$k\in{\mathbb Z}$ note $\overline k=k+d{\mathbb Z}\in{\mathbb Z}_d$.
The nonassociative algebra $(\derA,[.,.]_{\sigma})$ is  ${\mathbb Z}_d$-graded:
$\derA=\oplus_{\overline k\in{\mathbb Z}_d}{\mathcal D}_{\overline
k}$, with ${\mathcal D}_{\overline k}= {\mathbb C}[T^{\pm1}]d_k$ for
any $k\in\overline k$.
\end{stheoreme}
\begin{demo}
The case $s=1$ is straightforward. So assume $s\neq 1$.
Then note that  $T=qt^d$ if $s>1$, and $T=qt^{-d}$ if $s<1$, with $d\geq 1$.
So ${\mathbb C}[t^{\pm 1}]=\oplus_{i=0}^{d-1}t^i{\mathbb C}[T^{\pm1}]$ as vector spaces.
Now the direct sum in the Theorem follows from this and from the fact that ${\rm Ann}\Delta=\{0\}$.
The grading results directly from formulas (\ref{bracketsgt1}) and
from  the fact that $t^{n+m}=t^{n+m-d}q^{-1}T$ if $s>1$, and $t^{n+m}=t^{n+m-d}qT^{-1}$ if $s<1$.
\end{demo}

%\bigskip

\subsection{The bracket ``modulo inner derivations"}\label{parmodulo}
Motivated by the previous results we are interested next in
obtaining more detailed description of what relations for the
bracket in \cite[Theorem 8]{HLS} (Theorem \ref{th:qWittnonlinear})
become modulo inner $\sigma$-derivations.

\noindent{\bf Notation.}
We will use the following notation for congruence of two $\sigma$-derivations
modulo inner $\sigma$-derivations: $\forall \ \wt\Delta, \widehat\Delta\in\derA$,
the expression $\wt\Delta\equiv\widehat\Delta $ means that $\wt\Delta-\widehat\Delta\in\innA$.

\medskip
%\showhyphens{symmetric}
Since $[\innA,\innA]_{\sigma}\subseteq\innA$, and the bracket is skew-sym\-met\-ric, we shall only compute modulo $\innA$ the brackets of 2 types:
\begin{itemize}
\item $[d_n,d_m]_{\sigma}$, with $0\leq n<m\leq d-1$;
\item $[d_n,gd_m]_{\sigma}$, with $m\in{\mathbb Z}$ and $0\leq n\leq d-1$.
\end{itemize}

Note that since $T=qt^{s-1}$ is a unit in ${\mathcal A}$, and for $s<1$ we have $1-T=-T(1-T^{-1})=-Tg$, so for all $s\in{\mathbb Z}$ we can write $\innA=(1-T)\derA$.

\begin{slemme} \label{anvers}
For all $n\in{\mathbb N}, \ m\in{\mathbb Z}$,  one has
\begin{enumerate}
\item $[d_n,d_m]_{\sigma}\equiv (n-m)d_{n+m}$;
\item if $s\geq1$ then  $d_m\equiv q^{-1}d_{m-d}$;\\
if $s<1$ then $d_m\equiv qd_{m-d}$.
\end{enumerate}
\end{slemme}
\begin{demo}
1. By Proposition \ref{331} we know that:\\
 $[d_n;d_m]_{\sigma}=-(T^n\sum_{k=0}^{m-n-1}T^k)d_{n+m}$.
 Since $\innA=(1-T)\calA\Delta$, we must compute the remainder of the Euclidian division of the polynomial in $T$ appearing in the right-hand side  by $1-T$.
Straightforward calculations show that
for all $p\in{\mathbb Z},\ r\in{\mathbb N},\ r\geq 1$,  we have in ${\mathbb C}[X^{\pm 1}]$
\begin{itemize}
\item $\ssum_{k=0}^rX^k=r+1-(1-X)(X^{r-1}+2X^{r-2}+\ldots+r)$;
\item if $p\geq 1$ then $X^p\ssum_{k=0}^rX^k=r+1-(1-X)(X^{p+r-1}+2X^{p+r-2}+\ldots+rX^p+(r+1)(X^{p-1}+\ldots+1))$;
\item if $p\leq -1$ then $X^p\ssum_{k=0}^rX^k=r+1+(1-X)((r+1)(X^{p}+\ldots+X^{-1})-X^{p+r-1}+2X^{P+r-2}+\ldots+rX^p)$.
\end{itemize}
It results from  these computations   that the remainder  is $(m-n)$.

2. Follows directly from the definition of $g$.
\end{demo}

\begin{sprop}\label{prop1}
 For all $0\leq n<m<d$ one has:
\begin{enumerate}
\item If $n+m<d$ then  $[d_n,d_m]_{\sigma}\equiv (n-m)d_{m+n}$. \\ In particular for $n=0$ one gets $[d_0,d_m]_{\sigma}=-md_m+\wt\Delta$, with $\wt\Delta\in\innA$, and for $m=1$ one has exactly $[d_0,d_1]_{\sigma}=-d_1$.
\item If $n+m\geq d$ then \begin{enumerate}
\item $[d_n,d_m]_{\sigma}\equiv {{(n-m)}\over q}d_{m+n-d}$ if $s\geq 1$;
\item $[d_n,d_m]_{\sigma}\equiv q(n-m)d_{m+n-d}$ if $s<1$.
\end{enumerate}\end{enumerate}
\end{sprop}
\begin{demo}
 If $n+m<d$ we are done thanks to Lemma \ref{anvers}, part 1.  One can easily check the case $n=0,m=1$. If $m+n\geq d$, we conclude using part 2 of Lemma \ref{anvers}.
\end{demo}

%\medskip

\begin{sremark}{\rm Note that in the first formula, there is no more $q$ appearing, just like in the classical Witt algebra. Unfortunately this does not induce such a formula on the quotient space $\derA/\innA$ because $\innA$ is not an ideal for the $[.,.]_{\sigma}$ bracket, as results from the next proposition.}
\end{sremark}

\begin{sprop}\label{prop2}
Assume $s\neq 1$.
 For all $0\leq n<d$, and all $m\in{\mathbb Z}$, set $p=\left[{{n+m}\over d}\right]$ the integral part of $(n+m)/d$. Then
 $[d_n,gd_m]_{\sigma}\equiv -dq^{-\epsilon p}d_{m+n-p d}$, with $\epsilon=\textrm{sign}(s-1)$.
\end{sprop}
\begin{demo}
 We compute $[d_n,gd_m]_{\sigma}=-[t^n\Delta,t^m\Delta]_{\sigma}+q[t^n\Delta,t^{m+d}\Delta]_{\sigma}$, with $0\leq n<d$ and $m\in{\mathbb Z}$.
 Case $s>1$: by Lemma \ref{anvers} we get $[d_n,d_m]_{\sigma}\equiv (n-m)d_{n+m}$ and $[d_n,d_{m+d}]_{\sigma}\equiv -(m+d-n)t^{m+d+n}\Delta$.
So
$$\begin{array}{c}[d_n,gd_m]_{\sigma}\equiv (n-m)d_{n+m} -q(n-m-d)d_{n+m+d}\\
\equiv(n-m)(1-qt^d)d_{m+n}-qdd_{n+m+d}\equiv -qdd_{n+m+d}\equiv -dd_{n+m}\Delta.\end{array}$$
 So now it depends on the value of $n+m$. Recall that $p$  is the only integer such that $dp\leq n+m<d(p+1)$. So $0\leq n+m-dp<d$, and thanks to Lemma \ref{anvers} part 2 we are done.

 Case $s<1$:
 $$\begin{array}{c} [t^n\Delta,t^mg\Delta]_{\sigma}=[t^n\Delta,t^m\Delta]_{\sigma}-q^{-1}[t^n\Delta,t^{m+d}\Delta]_{\sigma}\equiv\\ \equiv (n-m)t^{n+m}\Delta-q^{-1}(n-m-d)t^{n+m+d}\Delta\\
\equiv q^{-1}dt^{n+m+d}\Delta+(n-m)t^{n+m}(1-q^{-1}t^d)\Delta\equiv dt^{n+m}\Delta.\end{array}$$
Once again we conclude thanks to part 2 of Lemma \ref{anvers}.
\end{demo}

\begin{sremark} {\rm For $m=0$ we get that $[d_n,g\Delta]_{\sigma}\equiv -dd_n$,
from what we will deduce that $\{\widehat\Delta\in\derA\ |\ [\widehat\Delta,\innA]_{\sigma}\subseteq\innA\}=\innA$.}
\end{sremark}

\subsection{The spaces $S_1$ and $S^1$} \label{subsec:qWittS1}

Recall that $\derA=\calA\Delta$, $\innA=g\calA\Delta$ and $[\innA,\innA]_{\sigma}\subseteq \innA$.
The notations are the following: $\calA={\mathbb C}[t^{\pm 1}]$, $\sigma(t)=qt^s$, and  $g=1-\lambda t^d$, with $d=s-1$ and $\lambda=q$ if $s\geq 1$, and $d=1-s$ and $\lambda=q^{-1}$ if $s<1$.
We recall the following from Section \ref{secufd}.

\begin{sdefin}
$S^1={\rm Span}_{{\mathbb C}}[\innA,\innA]_{\sigma}$;\\
$\wt S_1=\{\wt\Delta\ |\ [\wt\Delta,S^1]_\sigma\subseteq \innA\}$;\\
$S_1=\{\wt\Delta\in\derA\ |\ [\wt\Delta;S^1]_{\sigma}\subseteq S^1\}$.
\end{sdefin}
Now we can describe these spaces in the case we are considering
here.

\begin{stheoreme} \label{thmssets}
\begin{enumerate}
\item If $s=1$ then $\derA=\innA=S^1=\wt S_1=S_1$.
\item If $s=0$ then $\derA={\mathbb C}\oplus\innA=\wt S_1=S_1$, and $S^1=0$.
\item If $s\neq0,1$ then $S^1\subset\innA$ is a strict inclusion, and $\wt S_1=\derA$.\\
Moreover, if $s\neq -1$ then $S_1=\innA$.
\end{enumerate}
\end{stheoreme}
\begin{demo}
1. This was already noted at the beginning of subsection
\ref{pardecompo}.

2. The decomposition of $\derA$ is Theorem \ref{thmdecompo}. For the rest, just note that in this case $g=1-q^{-1}t$ and $\sigma(g)=0$. So it follows from Lemma \ref{s1sup} that $S^1=0$. Then by definitions of these sets one gets $\wt S_1=S_1=\derA$.

3. We begin with $\wt S_1$. We give the proof in the case $s>1$, the case $s<0$ is treated exactly in the same way, while changing the formulas.  Note that $\sigma(g)=1-q^{d+1}t^{ds}$.
We prove now that for all $\wt\Delta\in\derA$ one has $[\wt\Delta,S^1]\subseteq \innA$. It is enough to show that $\forall n,m$, with $0\leq n<d$ and $m\in{\mathbb Z}$, one has $[t^n\Delta,\sigma(g)g\Delta]_{\sigma}\equiv 0$. But
$[t^n\Delta,\sigma(g)g\Delta]_{\sigma}=  [t^n\Delta,t^mg\Delta]_{\sigma}-q^{d+1}[t^n\Delta,t^{m+ds}g\Delta]_{\sigma}$. Thanks to Proposition \ref{prop2} the right hand side is congruent up to inner derivations to
$$dq^{-p}t^{m+n-pd}\Delta-q^{d+1}dq^{-p-s}t^{m+ds+n-pd-sd}\Delta=0,$$
 where $p=[{{m+n}\over d}]$.

 \smallskip

Now we come to $S_1$ for $s\neq0,1,-1$. Because $S^1\subset \innA$
and by definition of $S^1$ it is clear that $S_1\supseteq\innA$. For
the inverse inclusion, the following argument is valid for both
$s>1$ and $s<-1$, once we have noted that in both cases $g$ divides
$\sigma(g)$. We shall remark also that considering the Remark
\ref{rkk} for $q=t$ and $p=1$ we get $\Delta(t)g\sigma(g)\Delta\in
S^1$, with $\Delta(t)=t$ if $s>1$ and $\Delta(t)=-qt^{s}$ for
$s<-1$. Consider a $\sigma$-derivation $\tilde P(t)\Delta\in S_1$.
According to Theorem \ref{thmdecompo} one should write $\tilde
P(t)=\sum_{i=0}^{d-1}a_it^i+g(t)R(t)$, so that
$g(t)R(t)\Delta\in\innA$. We must show that
$P=\sum_{i=0}^{d-1}a_it^i=0$. We give details here for $s>1$, the
reader may check the case $s<-1$ in the same way. Because of Lemma
\ref{s1sup} and the preceding remark  we shall have
$[P(t)\Delta,tg\sigma(g)\Delta]_{\sigma}\in g\sigma(g)\calA\Delta$.
Since this bracket is equal to
$((\sigma(P)\Delta(tg\sigma(g))-\Delta(P)\sigma(tg\sigma(g)))\Delta$
and ${\textrm{Ann}}\Delta=0$, we have  $g\sigma(g)$ dividing $T=
\sigma(P)\Delta(tg\sigma(g))-\Delta(P)\sigma(tg\sigma(g))$. But
$\sigma(tg\sigma(g))=\sigma(t)\sigma(g)\sigma(\sigma(g))$ is a
multiple of $g\sigma(g)$ because $\sigma(g)$ is a multiple of $g$.
So $g\sigma(g)$ divides $\sigma(P)\Delta(tg\sigma(g))$. Now since
$\Delta$ is a $\sigma$-derivation we have
$$\Delta(tg\sigma(g))=qt^s\Delta(g\sigma(g))+\Delta(t)g\sigma(g),$$
so $g\sigma(g)$ divides $\sigma(P)qt^s\Delta(g\sigma(g))$, i.e. it
divides $\sigma(P)\Delta(g\sigma(g))$. Now, since
$\Delta(g\sigma(g))=\sigma(g)\Delta(\sigma(g))+\Delta(g)\sigma(g)$,
we get that $g$ divides the polynomial
$\sigma(P)(\Delta(\sigma(g))+\Delta(g))$.

Now we prove that $g$ and $Q=\Delta(\sigma(g))+\Delta(g)$ are
relatively prime, so by Gauss' Lemma $g$ must divide $\sigma(P)$. By
definition $\Delta=({\rm{id}}-\sigma)/g$, so
$\Delta(\sigma(g)+g)=(\sigma(g)-\sigma^2(g)+g-\sigma(g))/g=1-(\sigma^2(g)/g)$.
Once again it is conveniant to notice that since  $T=qt^d$, then
$g=1-T$, $\sigma(g)=1-T^s$ and $\sigma^2(g)=1-T^{s^2}$. So
$Q=-T\sum_0^{s^2-2} T^k$ and it is prime with $g$, because any root
$t_0\in{\mathbb C}$ of $g$ satisfies $T(t_0)=1$, so
$Q(t_0)=1-s^2\neq0$ (by our hypothesis on $s$).

Finally we are reduced to the hypothesis that $g$ divides $\sigma(P)$, with $P=\sum_{i=0}^{d-1}a_it^i$.
Then any root $t_0$ of $g$ must be a root of $P$.
Note that $g=1-\lambda t^d$ admits exactly $d$ distinct roots (the $d$-roots of $\lambda^{-1}$) in ${\mathbb C}$.
Let $t_0$ be one of these roots. Then $\sigma(P)(t_0)=\sum a_i(qt_0^s)^i=\sum a_it_0^i$ since $t_0^d=q^{-1}$
and $d=s-1$ if $s>1$, and $t_0^d=q$ and $d=1-s$ if $s<-1$. So $P(t_0)=0$, and $P$ admits $d$ distinct roots.
Because $P$ is of degree not higher than $d-1$ this implies $P=0$.
\end{demo}

%\medskip

\begin{sremark} {\rm \begin{enumerate}
\item Note that the computation of $S_1$ in the last case relies on the fact that $s^2\neq 1$,
and that is the reason why we could not treat the case $s=-1$;
\item Theorem \ref{thmssets} shows that for non-linearly $q$-deformed
Witt algebra, ``stabilizer-like'' sets $S_1, \wt S_1$  while baring
some information on relation between $\derA$ and $\innA$, do not
provide a chain of subalgebras between $\derA$ and $\innA$.
\end{enumerate}}
\end{sremark}

\section*{Acknowledgements}
{\small This research was initiated during a course delivered by the
second author at the Erasmus Intensive Program GAMAP: Geometric and
Algebraic Methods of Physics and Applications in University of
Antwerp in September 2005. Both authors wish to thank Prof. Fred Van
Oystaeyen for kind hospitality and fruitful discussions. The first
author thanks also the Centre for Mathematical Sciences at Lund
University for hospitality during his one week stay there in
February 2006.

The research was supported by the European Liegrits network, the
British Engineering and Physical Sciences Research Council
(EP/D034167/1), The Crafoord Foundation, The Swedish Foundation for
Cooperation in Research and Higher Education (STINT), The Royal
Physiographic Society in Lund and The Royal Swedish Academy of
Sciences.}

%\medskip

%The next step then is to define $S^2$, which could be the linear span of $[\innA,S^1]_{\sigma}$ or $[S^1,S^1]_{\sigma}$.


\begin{thebibliography}{99}

\bibitem{BG} Brown, K., Goodearl, K., Lectures on algebraic quantum groups.
Advanced Courses in Mathematics. CRM Barcelona. Birkh\"auser Verlag, Basel, 2002.

\bibitem{C} Cauchon, G., {\sl Effacement des d\'erivations
et spectres premiers des alg\`ebres quantiques}, J. Algebra {\bf 260} (2003), no. 2, 476--518.

\bibitem{ChaiIsKuLuk} Chaichian, M., Isaev, A. P., Kulish, P., Lukierski, J.,
\emph{$q$-deformed Jacobi identity, $q$-oscillators and $q$-deformed
infinite-dimensional algebras}, Phys. Lett. B \textbf{237} (1990),
no. 3-4, 401--406.

\bibitem{ChaiElinPop} Chaichian, M., Ellinas, D., Popowicz Z.,
\emph{Quantum conformal algebra with central extension}, Phys. Lett.
B \textbf{248} (1990), no. 1-2, 95--99.

\bibitem{ChaiPopPres} Chaichian, M., Popowicz, Z., Prenajder, P.,
\emph{q-Virasoro algebra and its relation to the q-deformed KdV
system}, Phys. Lett. B \textbf{249} (1990), no. 1, 63--65.

\bibitem{ChaiKuLukPopPresn} Chaichian, M., Kulish, P., Lukierski, J.,
Popowicz, Z., Prenajder, P., \emph{q-deformations of Virasoro
algebra and conformal dimensions}, Phys. Lett. B \textbf{262}
(1991), no. 1, 32--38.

\bibitem{CurtrZachos1} Curtright, T. L., Zachos, C. K.,
\emph{Deforming maps for quantum algebras}, Phys. Lett. B
\textbf{243} (1990), no. 3, 237--244.

\bibitem{DamKu} Damaskinsky, E. V., Kulish, P. P., \emph{Deformed oscillators and their applications}
(Russian),  Zap.\ Nauch.\ Semin.\ LOMI \textbf{189} 1991, 37--74;
Engl. transl. in  J. Soviet Math.  \textbf{62} no.~5, (1992),
2963--2986.

\bibitem{DelEtinQuantFieldStrings} Deligne, P., Etingof, P., Freed, D.S., Jeffrey, L.C., Kazhdan, D.,
\mbox{Morgan, J.W.,} Morrison, D.R., Witten, E., (Eds) \emph{Quantum
Fields and Strings: A Course for Mathematicians}, 2 vol., Amer.
Math. Soc., 1999.

\bibitem{DiFranMathiSen} Di Francesco, P., Mathieu, P., S\'en\'echal, D., \emph{Conformal
Field Theory}, Springer Verlag, 1997, 890 pp.

\bibitem{FLM} Frenkel, I., Lepowsky, J., Meurman, A., \emph{Vertex Operator Algebras and the Monster},
Academic Press, 1988, 508 pp.

\bibitem{Fuchs1} Fuchs, J., \emph{Affine Lie Algebras and Quantum Groups}, Cambridge University Press, 1992, 433 pp.

\bibitem{Fuchs2} Fuchs, J., \emph{Lectures on Conformal Field Theory and Kac-Moody Algebras},
Springer Lecture Notes in Physics 498, (1997), 1--54.

\bibitem{G} Goodearl, K.R., {\sl Prime ideals in Skew Polynomial Rings and Quantized Weyl Algebras},
J. Algebra {\bf 150} (1992), no. 2, 324--377.

\bibitem{HLS} Hartwig, J.T., Larsson, D., Silvestrov, S.D., {\sl Deformations of Lie algebras using $\sigma$-derivations},
J. Algebra {\bf 295} (2006), no. 2, 314--361.

\bibitem{HelSil-book}{\rm Hellstr{\"o}m, L., Silvestrov, S.D.}, \emph{Commuting
Elements in $q$-Deformed Heisenberg Algebras}, World Scientific,
Singapore, 2000, 256 pp. (ISBN: 981-02-4403-7).

\bibitem{Hu}  Hu, N. \emph{$q$-Witt algebras, $q$-Lie algebras, $q$-holomorph structure and representations},  Algebra Colloq.  {\bf 6}  (1999),  no. 1, 51--70.

\bibitem{Jakob}  Jakobsen, H.P., \emph{Matrix Chain Models and
   their $q$-deformations}, Preprint Mittag-Leffler Institute, Report
   No 23, 2003/2004, ISSN 1103-467X, ISRN IML-R- -23-03/04-SE.

\bibitem{JakobLee}  Jakobsen, H.P.,  Lee, H.C.-W.,
   \emph{Matrix Chain Models and Kac-Moody algebras},
in ``Kac-Moody Lie algebras and related topics",
  Contemp. Math. 343, Amer. Math. Soc., Providence, RI,
   (2004), 147--165.

\bibitem{K} Kassel, C., {\sl Cyclic homology of differential operators, the Virasoro algebra and a $q$-analogue},
Comm. Math. Phys.  {\bf 146}  (1992),  no. 2, 343--356.

\bibitem{Kabook} Kassel, C., \emph{Quantum Groups}, Graduate Texts in Mathematics, 155. Springer-Verlag, New York, 1995, 531 pp.

\bibitem{LS1} Larsson, D., Silvestrov, S.D.,
\emph{Quasi-hom-Lie algebras, Central Extensions and 2-cocycle-like
identities}, J. Algebra \textbf{288} (2005), 321--344.

\bibitem{LS2} Larsson, D., Silvestrov, S. D. , \emph{Quasi-Lie algebras}, in ``Noncommutative
Geometry and Representation Theory in Mathematical Physics",
Contemp. Math., 391, Amer. Math. Soc., Providence, RI, 2005,
241-248.

\bibitem{LS3} Larsson, D., Silvestrov, S.D.,
\emph{Quasi-deformations of $sl_2(\mathbb{F})$ using twisted derivations},
Preprint, available at \\
\texttt{http://arxiv.org/abs/math.RA/0506172}.

\bibitem{LeBruyn} Le Bruyn, L., \emph{Conformal $sl_2$ Enveloping Algebras}, Comm. Alg. \textbf{23} no. 4 (1995), 1325--1362.

\bibitem{LeBruynSmith} Le Bruyn, L., Smith, S.P., \emph{Homogenized $sl_2$}, Proc. AMS 118 (1993), 725--730.

\bibitem{LeBruynSmithvdBergh} Le Bruyn, L., Smith, S.P., Van den Bergh, M.,
\emph{Central Extensions of Three-Dimensional Artin--Schelter Regular Algebras}, Math. Zeit. \textbf{222} no. 2 (1996), 171--212.

\bibitem{LeBruynvdBergh} Le Bruyn, L., Van den Bergh, M.,
\emph{On Quantum Spaces of Lie Algebras}, Preprint, UIA, available at \texttt{http://www.math.ua.ac.be/$\sim$lebruyn/}.

\bibitem{MCR} McConnell, J.C., Robson, J.C., Noncommutative Noetherian rings. Revised edition.
Graduate Studies in Mathematics, 30. American Mathematical Society, Providence, RI, 2001.

\bibitem{RS} Redondo, M.-J., Solotar, A., {\sl $\alpha$-derivations},  Canad. Math. Bull.  {\bf 38}  (1995),  no. 4, 481--489.

\end{thebibliography}
\end{document}